\documentclass[11pt]{article}
\usepackage{amsmath,amsfonts,amsthm}
\usepackage{xypic}

\def\P{{{\mathbb P}^{\ell-1}}}

\def\BB{{\mathbb B}}
\def\CC{{\mathbb C}}

\def\KK{{\mathbb K}}

\def\cA{{\cal A}}
\def\cG{{\cal G}}

\def\cB{{\cal B}}

\def\cC{{\cal C}}

\def\cO{{\cal O}}

\def\cV{{\cal V}}

\def\cK{{\cal K}}

\newtheorem{theorem}{Theorem}

\newtheorem*{theo*}{Theorem}
\newtheorem{lemma}{Lemma}[section]
\newtheorem{corol}[lemma]{Corollary}
\newtheorem{prop}[lemma]{Proposition}

\theoremstyle{definition}
\newtheorem{exam}[lemma]{Example}
\newtheorem*{defi}{Definition}

\DeclareMathOperator{\tot}{tot}
\DeclareMathOperator{\cut}{cut}
\DeclareMathOperator{\depth}{depth}
\DeclareMathOperator{\Ind}{Ind}
\DeclareMathOperator{\Indhom}{Ind_{hom}}

\DeclareMathOperator{\IndGSV}{Ind_{GSV}}
\DeclareMathOperator{\Img}{Im}
\DeclareMathOperator{\coker}{coker}
\DeclareMathOperator{\rank}{rank}
\DeclareMathOperator{\Ann}{Ann}

\newcommand{\indV}{\Ind_{{\rm Hom},V}}

\begin{document}
\title{\bf An Algebraic Formula for the Index of a Vector Field on an Isolated Complete Intersection Singularity
\thanks{Supported by DFG, Program "Global Methods in Complex Geometry" (Eb 102/4-3).}}
\author{\bf H.-Ch.\,Graf von Bothmer, W.\,Ebeling and X.\,G\'omez-Mont}
\maketitle

\begin{abstract}
Let $(V,0)$ be a germ of a complete intersection variety in $\CC^{n+k}$, $n>0$,  having an isolated singularity at $0$  and $X$ be
the germ of a holomorphic vector field on $\CC^{n+k}$ tangent to $V$ and having on $V$ an isolated zero at $0$.  We show that in this case the homological index and the GSV-index coincide. 
In the case when the zero of $X$ is also isolated in the ambient space $\CC^{n+k}$ we give a
formula for the homological index in terms of local linear algebra.
\end{abstract}

\section*{Introduction}
An isolated singular point (zero) $p$ of a vector field on $\CC^n$ has an index. It can be defined as  the degree of the map $X/\| X \|$ from a small sphere around the point $p$ to the unit sphere. 
If the vector field is holomorphic, then the index can also be defined as the dimension of a certain algebra: If, in local coordinates centred at the point $p$,
$$X= \sum_{i=1}^n X_i \frac{\partial }{\partial x_i},$$
then the index is equal to the dimension of the complex vector space $\cO_{\CC^n,0}/(X_1, \ldots , X_n)$, where $\cO_{\CC^n,0}$ is the ring of germs of holomorphic functions of $n$ variables and $(X_1, \ldots , X_n)$ is the ideal generated by the components of the vector field $X$.

Now let $X$ be a vector field tangent to the germ $(V,0)$ of a complex analytic variety of pure dimension $n$ with an isolated singularity at $0$ such that $X$ has an isolated singular point at 0 as well. Then one can try to generalize the two notions of the index mentioned above. If $V$ is a complete intersection, then one can  still define an index of $X$ at  $0$ as the degree of a certain map. This is done in \cite{GSV, SS} and this index is called the GSV-index $\IndGSV (X;V,0)$. If, more generally,  $V$ is a complex analytic variety of pure dimension $n$, but $X$ is a holomorphic vector field, then the third author proposed a generalization of the algebraic index, which is called the {\em homological index}  \cite{G}. This is defined as follows: Consider the sheaves $\Omega_V^j$ of germs of differential $j$-forms on $V$, $j=0, \ldots , n$. Contraction by the vector field $X$ defines a complex
\[
	\Omega_{V} \colon 
	0 \xleftarrow{}
	\cO_V \xleftarrow{X}
	\Omega^1_V \xleftarrow{X}
	\Omega^2_V \xleftarrow{X}
	\dots \xleftarrow{X}
	\Omega^n_V \xleftarrow{}
	0
\]
Since both $(V,0)$ and $X$ have isolated singularities, this complex has finite dimensional homology groups. The {\em homological index} $\Indhom (X; V,0)$ of the vector field $X$ at 0  is defined as the Euler characteristic of this complex
$$\Indhom (X;V,0) := \chi(\Omega_V) = \sum_{i=0}^n (-1)^i \dim H_i(\Omega_V).$$
In \cite{G} it was shown that the homological index and the GSV-index differ only by a constant depending on the germ of the variety $(V,0)$ but not on the vector field $X$. In the case when $V$ is a hypersurface in $\CC^{n+1}$  it was shown in that paper that this constant is equal to zero. Moreover, in this case an algebraic formula for the homological index under the additional hypothesis that $X$ has an isolated singularity in $\CC^{n+1}$ was given.


In this paper we consider the case when  $(V,0)$ is the germ of a complete intersection variety having an isolated singularity at $0$ defined by the vanishing of the germs of holomorphic functions $f_1,\ldots,f_k$ in $\CC^{n+k}$, $n>0$, and $X$
the  germ of a holomorphic vector field on $\CC^{n+k}$ tangent to $V$ and having on $V$ an isolated zero at $0$. We show that in this case the homological index and the GSV-index  coincide (Theorem \ref{tConstant0}). 
In the case when the zero of $X$ is also isolated in the ambient space $\CC^{n+k}$ we give a formula for the homological index in terms of local linear algebra (Corollary \ref{cSingularity}).
When $V$ is a hypersurface we recover the formula of \cite{G}. 

Our method of proof is as follows. The tangency condition can be expressed by an anticommutative square of finite free $\cO$-modules. Given such a square we construct a double complex which we call the Gobelin\footnote{A Gobelin is a richly embroidered French wall tapestry.}. It turns out that the Gobelin is weaved from Koszul-complexes and complexes introduced by Buchsbaum and Eisenbud. Using their results  and
the first spectral sequence of the Gobelin we prove that one can cut the Gobelin to obtain a finite free resolution of the complex $\Omega_V$. It follows that the Euler characteristic of the cut Gobelin is equal to the homological index of $V$. We then construct a deformation $V_\lambda$ of $V$ to its Milnor fiber and a family of vector fields $X_\lambda$ tangent to $V_\lambda$. Applying the Gobelin construction to this family yields a situation in which we can apply the results of \cite{GG02} and \cite{G} and hence conclude that the Euler characteristics of the family of cut Gobelins is independent of $\lambda$. 
It follows that the GSV-Index and the homological index agree. 

The second spectral sequence of the Gobelin provides formulae for the homology groups of the Gobelin in terms of local linear algebra.

The paper is organized as follows. We first define the Gobelin double complex $\cG$ and show in Theorem~\ref{tCut} that the vertical complexes of the Gobelin are resolutions up to a certain column.
In Section~\ref{sConstant0} we apply this to the situation that $X$ is a holomorphic vector field tangent to a germ of  an $n$-dimensional complete intersection singularity $(V,0)$ in $\CC^N$ with $X$ and $V$ both having an isolated singularity at $0 \in \CC^N$. We show that the homological index of $X$ at 0  is the Euler characteristic of the total complex of the subcomplex of the Gobelin consisting of the first $n+1$ columns of $\cG$ (Proposition~\ref{pCutIndex}). We derive from this that the homological index coincides with the GSV-index. In order to compute the homological index we show in Section~\ref{ssmallGobelin} that the total complex of the Gobelin is quasi-isomorphic to the total complex of a simpler double complex which we call the small Gobelin. Here we need the fact that $X$ has an isolated singularity in the ambient space $\CC^N$. In Section~\ref{sFormulae} we derive from this algebraic formulae for the homological index. The main formula is contained in Corollary~\ref{cSingularity}. 

We wish to thank the referee for useful comments, in particular for pointing out an error in our previous calculation of Example~Ê\ref{Ex2}.

\section{The Gobelin Double Complex} \label{sGobelin}

Throughout this paper we fix a field $\KK$ and consider a local Noetherian $\KK$-algebra $\cO$. All tensor products will be over $\KK$, unless otherwise specified. If $\varphi$ is a matrix with entries in $\cO$, we denote by $I_\ell(\varphi)$ the ideal of its $\ell \times \ell$ minors.
A complex has {\em length} $r$ if it contains $r+1$ non-zero elements. A complex is called a {\em resolution} if its first non zero term is in degree 0 and it has zero homology in positive degrees.

In this section we will develop the technical tools needed to prove our results. 

\subsection{Construction of the Gobelin}

Consider a matrix identity over $\cO$

\begin{equation}
\begin{pmatrix}\varphi_{11} & \cdots& \varphi_{1N}\cr
\vdots &\ddots& \vdots \cr
\varphi_{\ell 1} & \cdots& \varphi_{\ell N}\cr
\end{pmatrix}
\begin{pmatrix}X_1\cr
 \vdots \cr
X_N\cr
\end{pmatrix}
=
\begin{pmatrix}c_{11} & \cdots& c_{1k}\cr
\vdots &\ddots& \vdots \cr
c_{\ell 1} & \cdots& c_{\ell k}\cr
\end{pmatrix}
\begin{pmatrix}f_1\cr
\vdots \cr
f_k\cr
\end{pmatrix}
\label{(0.1)}
\end{equation}
and write this equation  as
\begin{equation}
\varphi X = c f.
\label{E0}
\end{equation}
Let
$F$, $G$ and $H$ be finite dimensional $\KK$-vector spaces of dimensions $N$, $\ell$, and $k$ respectively. Then the equation (\ref{E0})
gives rise to the anticommutative square
\begin{equation}
\begin{matrix}
\cO &
{\overset{-f}{\longrightarrow}}&
H \otimes\cO\cr
X \downarrow \ \ &&\downarrow
 c\cr
F \otimes \cO&
{\overset{\longrightarrow}{ \varphi}}
&
G \otimes \cO \cr
\end{matrix}
\label{E1}
\end{equation}
Let $\P$ denote the projective space ${\rm Proj}(G)$ and $\cO_{\P}(1)$ the sheaf of hyperplane sections on $\P$. Let $s_1,\ldots,s_\ell$ be a basis of its global sections, $s:=(s_1,\ldots,s_\ell)$,
$\widetilde \cO :=\cO \otimes\cO_{\P}$ and
$\widetilde \cO(m) :=\cO \otimes\cO_{\P}(1)^{\otimes m}$.
We tensor the diagram (\ref{E1}) with the sheaf $\cO_\P$
and continue at the right bottom of the square
with the tensor product of the natural morphism
$$s\cdot\  :G\otimes \cO_\P \longrightarrow \cO_\P(1)$$
with $\cO$ to obtain the following anticommutative square
of $\widetilde\cO$-sheaves on $\P$:
\begin{equation}
\begin{matrix}
\widetilde\cO &
{\overset{-f}{\longrightarrow}}&
H \otimes\widetilde\cO\cr
X \downarrow \ \ &&\downarrow
s \cdot c\cr
F \otimes \widetilde\cO&
{\overset{\longrightarrow}{s \cdot \varphi}}
&
\widetilde\cO(1) \cr
\end{matrix}
\label{(0.1a)}
\end{equation}
Since going around the square gives a $1  \times 1$ matrix, we may transpose the upper part of the square and obtain the anticommutative square
\newcommand{\smallXc}{\left(\begin{smallmatrix}X\\ (s\cdot c)^t\end{smallmatrix}\right)} 
\begin{equation}
\begin{matrix}
\widetilde\cO &
\xrightarrow{(s \cdot c)^t}&
H^* \otimes\widetilde\cO(1)\cr
X \downarrow \ \ &&\downarrow
-f^t\cr
F \otimes \widetilde\cO&
{\overset{\longrightarrow}{s \cdot \varphi}}
&
\widetilde\cO(1) \cr
\end{matrix}.
\label{(0.2)}
\end{equation}
The total complex associated to this square is the syzygy
\begin{equation}
\widetilde\cO
\xrightarrow{\smallXc}
[F \otimes \widetilde\cO] \oplus
[H^* \otimes\widetilde\cO(1)] 
\xrightarrow{(s \cdot \varphi , -f^t)}
\widetilde\cO(1)
\label{(0.3)}
\end{equation}
Using this
and the notation
\begin{equation}
\cV:= [F \otimes \widetilde\cO] \oplus
[H^* \otimes\widetilde\cO(1)],
\quad \cV(m):=\cV \otimes_{\widetilde \cO} \widetilde \cO(m)
\label{(0.3a)}
\end{equation}
we consider upward vertical Koszul complexes  
$\cK_{(s \cdot \varphi ,-f^t)} \otimes \cO(-k+i-1)$ of length $N+k$,
\begin{equation}
0 \leftarrow  \cO(-k+i)  \leftarrow
 \cV(-k+i-1)  \leftarrow
 \Lambda^2\cV(-k+i-2) \leftarrow \, \cdots \, \leftarrow
 \Lambda^{N+k} \cV(-k+i-N-k)
\leftarrow 0
\end{equation}
whose arrows are given by contractions with
$(s \cdot \varphi ,-f^t)$ and 
rightward horizontal Koszul complexes 
$\cK_{\smallXc} 
\otimes \cO(-k+j)$ of length $N+k$,
\begin{equation}
0 \rightarrow  \cO(-k+j)  \rightarrow
\cV(-k+j)  \rightarrow
 \Lambda^2\cV(-k+j) \rightarrow \, \cdots \,
\rightarrow
 \Lambda^{N+k}\cV(-k+j)
\rightarrow 0
\end{equation}
with arrows given by exterior products with ${\smallXc}$.
These complexes we weave into a bi-infinite double complex of sheaves
$\widetilde \cG$
on $\P$ (see Fig.~\ref{fdoubleK} for an example). The anticommutativity of this double complex 
follows from the anticommutativity of exterior product and contraction
in the exterior algebra and the fact that the syzygy (\ref{(0.3)}) of this
square is a complex.

\newcommand{\bolduparrow}{\ \ \boldsymbol \uparrow}
\newcommand{\boldlongrightarrow}{\boldsymbol \longrightarrow}
\newcommand{\bcV}{\boldsymbol \cV}
\newcommand{\bcO}{\boldsymbol \cO}

\begin{figure}[h]
 {\scriptsize
$$
\begin{array}{*{18}{c@{\,}}}

&
&&
&&
& \uparrow &
& \uparrow &
& \uparrow &
& \uparrow &
& \uparrow &
& \uparrow &
\\

&
&&
&{\mathbf 0}  & \longrightarrow
& \widetilde \cO(2)  & \longrightarrow
& \cV(2)  & \longrightarrow
&\Lambda^2 \cV(2)  & \longrightarrow
&\Lambda^3 \cV(2)  & \longrightarrow
&\Lambda^4 \cV(2)  & \longrightarrow
&\Lambda^5 \cV(2)  & \longrightarrow
\\
&
&&
& \bolduparrow &
& \uparrow &
& \uparrow & 
& \uparrow &
& \uparrow &
& \uparrow &
\\

&&
0        & \longrightarrow
&\mathbf{\widetilde \bcO(1)} & \longrightarrow 
&\cV(1)  & \longrightarrow
&\Lambda^2\cV(1)  & \longrightarrow
&\Lambda^3 \cV(1)  & \longrightarrow
&\Lambda^4 \cV(1)  & \longrightarrow
&\Lambda^5 \cV(1)   & \longrightarrow
&\Lambda^6 \cV(1)   & \longrightarrow
\\

&
& \uparrow &
& \bolduparrow &
& \uparrow & 
& \uparrow &
& \uparrow &
& \uparrow &
& \uparrow &
& \uparrow &

\\
{\bf 0          }& {\boldlongrightarrow}
&{\bf \widetilde \bcO}  & {\boldlongrightarrow}
&{\bf\bcV}  & {\boldlongrightarrow}
& {\bf\Lambda^2 \bcV}  & {\boldlongrightarrow}
&{\bf\Lambda^3\bcV}  & {\boldlongrightarrow}
&{\bf\Lambda^4\bcV}  & {\boldlongrightarrow}
&{\bf\Lambda^5\bcV}  & {\boldlongrightarrow}
&{\bf\Lambda^6\bcV}   & {\boldlongrightarrow} 
& \mathbf{0}\\

\uparrow &
& \uparrow &
& \bolduparrow & 
& \uparrow &
& \uparrow &
& \uparrow &
& \uparrow &
& \uparrow &
\\

\widetilde \cO(-1)  & \longrightarrow
&\cV(-1)  & \longrightarrow
& {\bf\Lambda^2 \bcV(-1) } & \longrightarrow
&\Lambda^3\cV(-1)  & \longrightarrow
&\Lambda^4\cV(-1)  & \longrightarrow
&\Lambda^5\cV(-1)  & \longrightarrow
&\Lambda^6\cV(-1)   & \longrightarrow
& 0
&\\

\uparrow &
& \uparrow &
& \bolduparrow & 
& \uparrow &
& \uparrow &
& \uparrow &
&\uparrow
\\

\cV(-2)  & \longrightarrow
&\Lambda^2\cV(-2)  & \longrightarrow
&{\bf\Lambda^3 \bcV(-2)}  & \longrightarrow
&\Lambda^4 \cV(-2)  & \longrightarrow
&\Lambda^5 \cV(-2)   & \longrightarrow
&\Lambda^6 \cV(-2)   & \longrightarrow
& 0&
\\

 \uparrow &
& \uparrow &
& \bolduparrow & 
& \uparrow &
& \uparrow &
& \uparrow &
&

\\

\Lambda^2 \cV(-3)  & \longrightarrow
&\Lambda^3\cV(-3)  & \longrightarrow
&{\bf\Lambda^4\bcV(-3) } & \longrightarrow
&\Lambda^5\cV(-3)  & \longrightarrow
&\Lambda^6\cV(-3)   & \longrightarrow
&0& \\

 \uparrow & 
& \uparrow &
& \bolduparrow &
& \uparrow &
& \uparrow &

\\
 \Lambda^3 \cV(-4)  & \longrightarrow
&\Lambda^4 \cV(-4)  & \longrightarrow
&{\bf\Lambda^5\bcV(-4) } & \longrightarrow
&\Lambda^6\cV(-4)   & \longrightarrow
& 0&
\\

\uparrow &
& \uparrow &
& \bolduparrow &
& \uparrow

\\

\Lambda^4 \cV(-5)  & \longrightarrow
&\Lambda^5\cV(-5)  & \longrightarrow
&{\bf\Lambda^6\bcV(-5)}   & \longrightarrow
&0&
\\

 \uparrow &
& \uparrow &
& \bolduparrow & 
&

\\
 \Lambda^5\cV(-6)  & \longrightarrow
&\Lambda^6\cV(-6)   & \longrightarrow
&{\bf 0}&
\\

 \uparrow &
& \uparrow &
&

\end{array}
$$
}

\begin{center}
\parbox{12cm}{
 \caption{Part of the double complex $\widetilde \cG$ for $N=4$, $k=2$ and $l=2$. The generating Koszul complexes of (\ref{(0.3)}) are typed in a darker tone.} \label{fdoubleK}
 } \end{center}
\end{figure}


\begin{defi}
We define the {\em Gobelin double complex} 
of the anticommutative square (\ref{E1}) $$\cG := \Gamma(\widetilde \cG)^*$$  as the
dual of the global section double complex of $\widetilde \cG$.
The terms of the Gobelin are $\cO$-modules
$$ \cG_{i,j} = H^0(\P, \Lambda^{k+i-j}\cV(-k+j))^* = \bigoplus_{r+s=k+i-j}D_{-k+j+r}G^* \otimes\Lambda^rH \otimes
\Lambda^s F^*\otimes  \cO$$
where $D_mG^*:=H^0(\P,\cO_{\P}(m))^*$ is the homogeneous component of the divided power algebra
of $\KK[x_1,\ldots,x_\ell]$ of degree $m$, and the connecting maps are constructed using the matrices
$f$ and $\varphi$ for the vertical strands and $X$ and $c$ for the horizontal ones (see Fig.~\ref{fGobelin} for an example).
\end{defi}



\renewcommand{\textfraction}{0.1}
\begin{figure}[t]
 {\tiny
$$
\begin{array}{*{18}{c@{\,}}}

&
&&
&&
&&
&&
& \downarrow &
& \downarrow &
\\
\\

&&
&&
&&
&&0 & \longleftarrow
& D_4G^* \otimes \cO & \longleftarrow
&\begin{matrix}
D_5G^*\otimes H\otimes  \cO \\
\oplus\\
D_4G^*\otimes F^* \otimes \cO
\end{matrix}

\\ \\

&
&&
&&
&&
& \downarrow &
& \downarrow &
& \downarrow &

\\ \\

&&&&&& 0 & \longleftarrow
&D_3G^* \otimes \cO & \longleftarrow
&\begin{matrix}
D_4G^*\otimes H\otimes  \cO \\
\oplus\\
D_3G^*\otimes F^* \otimes \cO
\end{matrix}
 & \longleftarrow
& \begin{matrix}
D_5G^*\otimes\Lambda^2H \otimes \cO \\
\oplus\\
D_4G^*\otimes H\otimes  F^* \otimes \cO \\
\oplus\\
D_3G^*\otimes\Lambda^2 F^* \otimes \cO
\end{matrix}
&

\\ \\

&
&&
&  &
& \downarrow &
& \downarrow &
& \downarrow &
& \downarrow &

\\ \\

&&&&
 0 & \longleftarrow
&D_2G^*\otimes\cO & \longleftarrow
&\begin{matrix}
D_3G^*\otimes H  \otimes \cO \\
\oplus\\
D_2G^*\otimes F^* \otimes \cO
\end{matrix}
 & \longleftarrow
&\begin{matrix}
D_4G^*\otimes\Lambda^2H  \otimes \cO \\
\oplus\\
D_3G^*\otimes H\otimes  F^* \otimes \cO \\
\oplus\\
D_2G^*\otimes\Lambda^2 F^* \otimes \cO
\end{matrix}
 & \longleftarrow
&\begin{matrix}
D_4G^*\otimes\Lambda^2H \otimes  F^* \otimes \cO \\
\oplus\\
D_3G^*\otimes H\otimes \Lambda^2 F^* \otimes \cO \\
\oplus\\
D_2G^*\otimes\Lambda^3 F^* \otimes \cO
\end{matrix}

\\ \\

&&&
& \downarrow &
& \downarrow &
& \downarrow &
& \downarrow &
& \downarrow &

\\ \\

&&
0          & \longleftarrow
& G^* \otimes \cO & \longleftarrow
&\begin{matrix}
D_2G^*\otimes H \otimes \cO \\
\oplus\\
G^*\otimes  F^* \otimes \cO
\end{matrix}
 & \longleftarrow
&\begin{matrix}
D_3G^*\otimes\Lambda^2H \otimes \cO \\
\oplus\\
D_2G^*\otimes H\otimes F^* \otimes \cO \\
\oplus\\
G^*\otimes\Lambda^2 F^* \otimes \cO
\end{matrix}
 & \longleftarrow
&\begin{matrix}
D_3G^*\otimes\Lambda^2H \otimes F^* \otimes \cO \\
\oplus\\
D_2G^*\otimes H\otimes \Lambda^2 F^* \otimes \cO \\
\oplus\\
G^*\otimes\Lambda^3 F^* \otimes \cO
\end{matrix}
 & \longleftarrow
& \begin{matrix}
D_3G^*\otimes\Lambda^2H \otimes \Lambda^2 F^* \otimes \cO \\
\oplus\\
D_2G^*\otimes H\otimes \Lambda^3 F^* \otimes \cO \\
\oplus\\
G^*\otimes\Lambda^4 F^* \otimes \cO
\end{matrix}

\\ \\

&&\downarrow  &
& \downarrow &
& \downarrow &
& \downarrow &
& \downarrow &
& \downarrow &

\\ \\

0          & \longleftarrow&
 \cO & \longleftarrow
&\begin{matrix}
G^*\otimes  H\otimes \cO \\
\oplus\\
 F^* \otimes \cO
\end{matrix}
 & \longleftarrow
&\begin{matrix}
D_2G^*\otimes\Lambda^2H \otimes \cO \\
\oplus\\
G^*\otimes H\otimes F^* \otimes \cO \\
\oplus\\
\Lambda^2 F^* \otimes \cO
\end{matrix}
 & \longleftarrow
&\begin{matrix}
D_2G^*\otimes\Lambda^2H \otimes F^* \otimes \cO \\
\oplus\\
G^*\otimes H\otimes \Lambda^2 F^* \otimes \cO \\
\oplus\\
\Lambda^3 F^* \otimes \cO
\end{matrix}
& \longleftarrow
&\begin{matrix}
D_2G^*\otimes\Lambda^2H \otimes \Lambda^2 F^* \otimes \cO \\
\oplus\\
G^*\otimes H\otimes \Lambda^3 F^* \otimes \cO \\
\oplus\\
\Lambda^4 F^* \otimes \cO
\end{matrix}
 & \longleftarrow
&\begin{matrix}
D_2G^*\otimes\Lambda^2H \otimes \Lambda^3 F^* \otimes \cO \\
\oplus\\
G^*\otimes H\otimes \Lambda^4 F^* \otimes \cO \\
\end{matrix}

\\ \\

 &&\downarrow &
& \downarrow &
& \downarrow &
& \downarrow &
& \downarrow &
& \downarrow &

\\ \\

0          & \longleftarrow&
 H\otimes \cO & \longleftarrow
&
\begin{matrix}
G^*\otimes\Lambda^2H\otimes \cO \\
\oplus\\
H\otimes F^* \otimes \cO
\end{matrix}
& \longleftarrow
&\begin{matrix}
G^*\otimes\Lambda^2H\otimes F^* \otimes \cO \\
\oplus\\
H\otimes \Lambda^2 F^* \otimes \cO
\end{matrix}
 & \longleftarrow
&\begin{matrix}
G^*\otimes\Lambda^2H\otimes \Lambda^2 F^* \otimes \cO \\
\oplus\\
H\otimes \Lambda^3 F^* \otimes \cO
\end{matrix}
 & \longleftarrow
&\begin{matrix}
G^*\otimes\Lambda^2H\otimes \Lambda^3 F^* \otimes \cO \\
\oplus\\
H\otimes \Lambda^4 F^* \otimes \cO
\end{matrix} & \longleftarrow
&
G^*\otimes\Lambda^2H\otimes \Lambda^4 F^* \otimes \cO \\

\\ \\

&& \downarrow &
& \downarrow &
& \downarrow &
& \downarrow &
& \downarrow &
&\downarrow &
&  &
&  &

\\ \\

0          & \longleftarrow&
 \Lambda^2H\otimes \cO & \longleftarrow
&\Lambda^2H\otimes F^* \otimes \cO & \longleftarrow
&\Lambda^2H\otimes \Lambda^2 F^* \otimes \cO & \longleftarrow
&\Lambda^2H\otimes \Lambda^3 F^* \otimes \cO & \longleftarrow
&\Lambda^2H\otimes \Lambda^4 F^* \otimes \cO & \longleftarrow
&0

&
\\ \\

 &&\downarrow &
& \downarrow &
& \downarrow &
& \downarrow &
&\downarrow &

\\ \\

&&0&
& 0&
& 0&
& 0&
& 0&
\\ \\
\end{array}
$$}
\caption{The lower left hand part of the Gobelin for $N=4$, $k=2$, $l=2$
beginning at $(0,0)$.} \label{fGobelin}
\end{figure}

\begin{lemma} \label{lem1}
Let $\cG$ be the Gobelin constructed from the identity
(\ref{(0.1)}) over the local $\KK$-algebra $\cO$. We have:

\begin{enumerate}
\item
The Gobelin is a double complex
of finite free $\cO$-modules and $\cG_{i,j}$ is non-zero only for
$j \geq0,\ i = \max\{0,j-k\},\ldots,N+j$.
\item
The $0^{th}$ row of the Gobelin is the complex
$$0          \longleftarrow
 \Lambda^kH \otimes \cO \longleftarrow
\Lambda^kH \otimes F^* \otimes \cO \longleftarrow
\Lambda^kH \otimes \Lambda^2  F^* \otimes \cO \longleftarrow
\cdots \longleftarrow
\Lambda^kH \otimes \Lambda^N  F^* \otimes \cO \longleftarrow
0
$$
with maps being  contractions with $X$, i.e.\ the tensor product of $\Lambda^k H$ with the Koszul complex $\cK_X$.

\end{enumerate}

\end{lemma}

\begin{proof} \hfill

\begin{enumerate}
\item
Since the Gobelin is obtained as the vector space of global sections of the
double complex of sheaves $\widetilde \cG$ and then taking duals, the anticommutativity
follows from the anticommutativity of $\tilde \cG$. Furthermore the summands of
$$ \cG_{i,j} =  \bigoplus_{r+s=k+i-j}D_{-k+j+r}G^* \otimes\Lambda^rH \otimes
\Lambda^s F^*\otimes  \cO$$
are nonzero if and only if $0 \le -k+j+r$ and $0 \le r \le k$ and $0 \le s \le N$.  

\item
This follows from substituting $j=0$ in the formula for $\cG_{i,j}$.

\end{enumerate}
\end{proof}


\subsection{The Vertical Complexes $\cA^m$ of the Gobelin}

\newcommand{\Kphi}{\cK_{s\cdot\varphi}}
\newcommand{\Kphimd}{(\Kphi^{m})^*}
\newcommand{\Kphim}{\Kphi^{m}}
\newcommand{\Kphinm}{\Kphi^{n-m}}

Denote by $\cA^m$ the complex obtained from the $m^{th}$ column of the Gobelin.
The complexes $\cA^m$ are the tensor product of two complexes, arising from the direct sum decomposition of the middle module in the syzygy (\ref{(0.3)}). One is the part of the Buchsbaum-Eisenbud
strand of the Koszul complex associated to the map
\begin{equation*}
\varphi:F\otimes \cO \longrightarrow G\otimes \cO
\end{equation*}
above the splicing map and the other is the
Koszul complex associated to the sequence $f_1,\ldots,f_k$.
In this subsection we apply the Buchsbaum-Eisenbud and Koszul Theorems
to describe the homologies of the vertical complexes in the Gobelin.

Denote by $\cK_f$ the Koszul complex obtained from the morphism $-f^t:H^* \otimes \cO \longrightarrow \cO$:

\begin{equation*}
0 \longleftarrow
  \cO  \longleftarrow
H^*  \otimes  \cO  \longleftarrow
\Lambda^2 H^* \otimes  \cO \longleftarrow
\cdots \longleftarrow
\Lambda^kH^* \otimes  \cO \longleftarrow
0
\end{equation*}
with morphisms contractions with $-f^t=(-f_1,\ldots,-f_k)$.
It's 0-homology group is
$$H_0(\cK_f) = \cO_V := \cO/(f_1,\ldots,f_k).$$
If $-f_1,\dots,-f_k$ is a regular $\cO$ sequence, then by Koszul's Theorem \cite[Theorem~17.4, p.~424]{E}
 $H_j(\cK_f)=0$ for $j>0$.
Let $\cK_f^*$ be its dual
$$
0 \longleftarrow
\Lambda^k H \otimes  \cO  \longleftarrow
\Lambda^{k-1} H  \otimes  \cO  \longleftarrow
\cdots \longleftarrow
H \otimes  \cO 
\longleftarrow \cO \longleftarrow
0
$$
with morphisms contractions with $-f^t=(-f_1,\ldots,-f_k)$.
It's 0-homology group is again
$$H_0(\cK_f^*) = \cO_V.$$
Next we consider the Koszul complex $\Kphi$ 
of $\widetilde \cO$-sheaves obtained from the morphism 
$$s\cdot\varphi: F \otimes  \widetilde \cO \longrightarrow \widetilde \cO(1)$$
on $\P$ and tensor it with $\widetilde \cO(m-1)$:
\begin{equation*}
0 \leftarrow  \widetilde \cO(m)  \leftarrow
 F \otimes \widetilde \cO(m-1)  \leftarrow
 \Lambda^2F \otimes \widetilde \cO(m-2) \leftarrow \cdots \leftarrow
 \Lambda^{N} F \otimes \widetilde \cO(m-N)
\leftarrow 0
\end{equation*}
As in \cite[A2.6.1, p.~591]{E} we denote by $\Kphim$ the complex of its global sections
\begin{equation*}
{
0 \leftarrow  S_{m}G \otimes \cO  \leftarrow
S_{m-1}G \otimes F \otimes \cO  \leftarrow
 S_{m-2}G \otimes \Lambda^2 F \otimes \cO \leftarrow \cdots \leftarrow
 S_{m-N}G \otimes \Lambda^N F \otimes \cO
\leftarrow 0}
\end{equation*}
where $S_mG=H^0(\P,\cO_\P(m))$,
 and  by $\Kphimd$
the dual complex of free $\cO$-modules
\begin{equation*}
{
0 \leftarrow  D_{m-N}G^* \otimes \Lambda^NF^* \otimes \cO  \leftarrow
  \cdots \leftarrow
 D_{m-1}G^* \otimes F^* \otimes \cO \leftarrow
 D_{m}G^*\otimes \cO
\leftarrow 0}
\end{equation*}
Note that for $m \le N$ the first nonzero term of $\Kphimd$ is $D_0 G^* \otimes \Lambda^m F^* \otimes \cO$. 
In this case 
the complexes $\Kphim$ and  $\Kphimd$ have length $m$. 
For $m \geq N$ the complex
$(\cK^m_{s\cdot\varphi})^*$ has length $N$.
For $m\leq N$ we have

$$H_0((\cK_{s\cdot\varphi}^m)^*) = \frac{\Lambda^mF^* \otimes \cO}
{ s\cdot \varphi (G^* \otimes \Lambda^{m-1}F^* \otimes \cO)} $$

$$ H^0(\cA^m) = \frac{\Lambda^mF^* \otimes \cO}
{ s\cdot \varphi (G^* \otimes \Lambda^{m-1}F^* \otimes \cO) + f( H^*\otimes \Lambda^mF^*
\otimes \cO)}$$

\begin{lemma} \label{lAgeneral}
For $m \ge 0$ we have
\begin{enumerate} 
\item
$\cA^m = \Kphimd \otimes_\cO\cK_f^*$


\item
If $f_1,\ldots,f_k$ is a regular $\cO$-sequence then 
$\cA^m$ is quasi-isomorphic to the complex
$\Kphimd\otimes_\cO \cO_V $ and
$H_j(\cA^m) = 0$ for $j>\min(N,m)$.

\end{enumerate}
\end{lemma} 
\begin{proof} \hfill

\begin{enumerate}
\item
Over $\P$, the $m$-th column $\cK_{(s \cdot \varphi ,-f^t)} \otimes \cO(m-1)$ of the double-complex $\widetilde \cG$ is a tensor product of two Koszul complexes, $[\Kphi \otimes {\widetilde \cO}(m-1)] \otimes_{\widetilde \cO} \cK_f$, due to the direct sum decomposition of $\cV$ in (\ref{(0.3a)}). Since the
second complex  is independent of the variables of $\P$, the tensor product can be taken over $\cO$.
Taking global sections and then dualizing we obtain $\cA^m = \Kphimd \otimes_\cO \cK^*_{f}$ as complexes.


\item
The homology of $\cA^m = \Kphimd \otimes \cK_f^*$ can be computed from the
double complex where we put in the horizontal axis the complex
$\Kphimd$ and on the vertical the complex $\cK_f^*$. If we compute the spectral sequence where we do first the vertical homology, we obtain by Koszul's Theorem that the only homology group is at $j=0$, where the homology complex is $\Kphimd \otimes \cO_V$. Since the spectral sequence degenerates, we obtain that the homology of this complex computes the homology of the double complex. This proves the first statement. Now the second statement is immediate from this, since the complex $\Kphimd \otimes\cO_V$ has length min$\{N,m\}$ which is shorter than the complex $\Kphimd \otimes \cK^*_f$, so the last homology groups of the larger complex
$\Kphimd \otimes \cK^*_f = \cA^m$ vanish.
\end{enumerate}
\end{proof}

\begin{lemma} \label{lACohenMac}
Assume that $n:=N-\ell \geq 0$ and that the depth of $I_\ell(\varphi) = n+1$,
the greatest possible value,
then we have:

\begin{enumerate}

\item
For $m=0,\ldots,n+1$ the complex $\Kphimd$ is a resolution of $H_0\bigl(\Kphimd\bigr)$.

\item
If $f_1,\ldots,f_k$ is a regular sequence in $\cO$ then $\cA^m$ is a resolution for
$m \leq n-k+1$.
For $m = n-k+2, \dots , n+1$ we have $H_j(\cA^m) = 0$ for $j>m-n+k-1$.
\item
If in addition $\cO$ is a local Cohen-Macaulay ring and $f_1,\ldots,f_r$ is an $\cO/I_\ell(\varphi)$-regular sequence then $\cA^m$ is a resolution for $m \leq n-k+r+1$
and for
$m=n-k+r+2, \dots , n+1$ we have $H_j(\cA^m) = 0$ for $j>m-n+k-r-1$.

\end{enumerate}
\end{lemma}

\begin{proof} \hfill

\begin{enumerate}
\item
For $m=0,\ldots,n$ we glue the complexes $\Kphinm$ on the left with $\Kphimd$ on the right using the splicing map
$\varepsilon: \Lambda^{m}F^*\cong \Lambda^{n-m+\ell}F \longrightarrow \Lambda^{n-m} F$
which is contraction by $\Lambda^\ell\varphi^t$:
$$0 \longleftarrow \Kphinm
{\overset{ \varepsilon}{\longleftarrow} }
\Kphimd \longleftarrow 0.$$
The complexes so obtained are called ${ \cC}^{n-m}$ in \cite[A2.6]{E}. There D. Eisenbud also defines
${ \cC}^{n-m} = \Kphinm$  for $m \le -1$ and
${ \cC}^{n-m} = \Kphimd$ for $m \ge n+1$.
The length of the complexes $\cC^{n-m}$ is $n+1$ for $m=-1,\ldots,n+1$.

The Buchsbaum-Eisenbud Theorem \cite[Theorem A2.10, p.594]{E} applied to $\varphi$
asserts that under our hypothesis the complex $\cC^{n-m}$ is a free resolution of $H_0(\cC^{n-m})$, for $m \le n+1$. 

If we cut the complex $\cC^{n-m}$ at the splicing map,
we obtain that for $m=0,\ldots,n+1$ the complex $\Kphimd$
is a free resolution of its $0$-homology module.

\item
For $m=0,\ldots,n+1$ consider the double complex
$ { \cC}^{n-m} \otimes {\cK}_f^*$  with
horizontal axis $i=0,\ldots,n+1$ and vertical axis $j=0,\ldots,k$.
Consider the spectral sequence where we first do vertical
homology. By Koszul's Theorem we only have non-zero terms for $j=0$, where the
homology is
${ \cC}^{n-m} \otimes { \cO}_V$.
Hence the
spectral sequence degenerates,
the total complex $\cC^{n-m} \otimes \cK_f^*$ is quasi-isomorphic to
$\cC^{n-m} \otimes \cO_V$
and the only non-zero
homology groups of the total complex are
in $j=0,\ldots,n+1$.

Now we do the other spectral sequence,  doing first the
horizontal homology. Again by the Buchsbaum-Eisenbud Theorem we obtain that the only non-vanishing terms are in
$i=0$ where we obtain the homology complex
$H_0({ \cC}^{n-m})  \otimes { \cK}_f^*$.
So again the spectral sequence
degenerates and
$\cC^{n-m} \otimes \cK_f^*$ is quasi-isomorphic to $H_0(\cC^{n-m}) \otimes \cK_f^*$.
Both spectral sequences together give $H_j(\cC^{n-m} \otimes \cK_f^*) = 0$ for $j>\hbox{min}\{k,n+1\}$.

If we cut the complex ${ \cC}^{n-m}$ at the splicing
map $\varepsilon$, the right hand side is $\Kphimd$.
The
double complexes  ${ \cC}^{n-m} \otimes {\cK}_f^*$
and $\cA^m = \Kphimd \otimes {\cK}_f^*$
coincide on the columns to the right of the splicing map. Since $\cC^{n-m}$ has a complex of length $n-m$ left of the splicing map, the $0$-th column of the cut double complex is the $(n-m+1)$-st column of the complex $\cC^{n-m}$.

Doing for both double complexes the vertical
homology first,
both spectral sequences degenerate, with 
$\cC^{n-m} \otimes \cO_V$ and $\Kphimd \otimes \cO_V$ respectively in the $0$th row.
Now, both of these complexes coincide to the right of the
splicing map, so that  we have
\[
	H_j(\cA^m) = H_j(\Kphimd \otimes \cO_V) = H_{j+n-m+1}(\cC^{n-m} \otimes \cO_V) \quad \text{for $j>0$.}
\]
Hence by the vanishing above, $\cA^m$ is a resolution of $H_0(\cA^m)$ if $n-m +2 > k$, i.e.\ $m \le n-k+1$.

For $m = n-k+2, \dots, n+1$ we still have $H_j(\cA^m) = 0$ for $n-m+1+j>k$, i.e.\ $j>m-n+k-1$.

\item Under the hypothesis that $\cO$ is Cohen-Macaulay,
we have by \cite[Corollary A2.13, p. 599]{E} that
$H_0(\cC^{n-m})$ is a maximal Cohen-Macaulay $\cO/I_\ell(\varphi)$-module for $m=0,\ldots,n+1$.  The assumption that $f_1,\ldots,f_r$ is a
$\cO/I_\ell(\varphi)$-regular sequence implies that
$f_1,\ldots,f_r$ is also a $H_0(\cC^{n-m})$-regular sequence by
\cite[Proposition 21.9, p.529]{E}.
This then means that $H_i\bigl(H_0(\cC^{n-m})\otimes\cK_f^*\bigr)=0$ for $i>k-r$ by Koszul's Theorem.
Repeating the argument at the end of Part 2 of this lemma,
we obtain that $\cA^m$ is a resolution for $n-m+2>k-r$ and that for
$m=n-k+r+2, \dots, n+1$ we have $H_j(\cA^m) = 0$ for $j>m-n+k-r-1$.
\end{enumerate}
\end{proof}



\begin{prop} \label{tA}
Let $\tot(\cG)$ be the total complex of the Gobelin $\cG$ 
constructed from the identity
$(\ref{(0.1)})$ over the local $\KK$-algebra $\cO$, and 
\[
H_0(\cA) \colon 
0 \longleftarrow
H_0(\cA^0) \longleftarrow
H_0(\cA^1) \longleftarrow
\cdots \longleftarrow
H_0(\cA^N) \longleftarrow
0
\]
the complex induced by taking vertical homology. Then one has the following statements:

\begin{enumerate}
\item
If $f_1,\ldots,f_k$ is an $\cO$ regular sequence,
$N \geq \ell$  and the depth of
$I_\ell(\varphi) =N-\ell+1$, the greatest possible value,  then
$$H_i\bigl(\tot(\cG)\bigr) = H_i\bigl(H_0(\cA)\bigr) \mbox{ for }i \le N-\ell-k+1.$$
\item
If in addition $\cO$ is Cohen-Macaulay  
and $f_1,\ldots,f_r$ is an $\cO / I_\ell(\varphi)$-regular sequence,  then
$$H_i(\tot(\cG)) = H_i(H_0(\cA)) \mbox{ for }i \le N-\ell-k+r+1.$$
\end{enumerate}
\end{prop}

\begin{proof}
We look at the spectral sequence, where we do vertical homology
first. By Lemma \ref{lACohenMac} the vertical strands $\cA^m$ have nonzero homology
only in step $0$ for $m \le N-\ell-k+1$ and for $m \le N-\ell-k+r+1$ using the stronger hypotheses.  
Since this spectral sequence converges to the homology of the total
Gobelin $\tot(\cG)$, this proves the claim.
\end{proof}

\begin{theorem} \label{tCut}
In the situation of Proposition 1.4 Part 2 for $i<N-\ell -k+r+1$
consider the finite double complexes $\cG^{\cut}_{\leq i} \subset \cG$
obtained by considering only the columns $\cA^m$ for $m=0,\ldots,i$.
Then $\tot(\cG^{\cut}_{\le i})$ is quasi-isomorphic to 
\[
H_0(\cA)_{\le i} \colon 
0 \longleftarrow
H_0(\cA^0) \longleftarrow
H_0(\cA^1) \longleftarrow
\cdots \longleftarrow
H_0(\cA^i) \longleftarrow
0
\]
\end{theorem}

\begin{proof}
By Lemma \ref{lACohenMac} {\em all} vertical strands $\cA^m$ of $\cG_{\leq i}^{\cut}$ are resolutions. 
\end{proof}


\section{Comparing Homological Index and GSV Index} \label{sConstant0}

Let $V$ be a germ of a complete intersection variety having an isolated singularity at $0$ defined by the vanishing of the germs of holomorphic functions $f_1,\ldots,f_k$ in $\CC^N$, $N>k$, and $X$
a germ of a holomorphic vector field having an isolated zero at $0$ in $\CC^{N}$
and tangent to $V$, i.e. $X(f) = c \cdot f$ where $c$ is the $k \times k$
matrix of cofactors. If we denote by $\varphi$ the Jacobi matrix
describing the differential $df$ of 
$f \colon \CC^N \to \CC^k$ we obtain the matrix equality $\varphi \cdot X = c \cdot f$.
Denote by $\cG$  the Gobelin constructed from this equality. Note that in this situation
$F^*  \otimes \cO= \Omega^1_{\CC^N}$ and $H = G$ are vector spaces of the same dimension $k$.

\begin{prop} \label{pCutIndex}
 In this situation set $n := \dim V = N-k$ and let $\cG^{\cut} := \cG^{\cut}_{\le n}$ be the subcomplex consisting of the first $n+1$
columns of $\cG$. Then 
$$\Indhom (X; V,0) = \chi(\tot(\cG^{\cut}))$$.
\end{prop}

\begin{proof}
We want to apply Theorem \ref{tCut}. For this let $\cO$ be the ring of convergent power series in $N$ variables. Since $V$ is a complete intersection $f_1,\dots,f_k$ is an $\cO$-regular sequence. 
Now $I_k(\varphi)$
describes the critical locus $C_f$ of $f$. Since $V$ has an isolated singularity 
the image of $C_f$ is a hypersurface in $\CC^k$. The critical locus $C_f$ has therefore codimension $N-k+1$ in $\CC^N$ and
$\depth I_k(\varphi)=N-k+1$ has the maximal possible value. Moreover $\cO$ is Cohen-Macaulay.
Since $V$ has an isolated singularity, the codimension of the singular locus in
the critical locus is $k-1$. After a holomorphic base change we can therefore assume
that $f_1, \dots, f_{k-1}$ is a regular sequence in $\cO/I_k(\varphi)$. So we can apply Theorem \ref{tCut} with $r=k-1$ and obtain that $\cG^{\cut}$ is quasi-isomorphic to 
$$
H_0(\cA)_{\le n} \colon 
0 \longleftarrow
H_0(\cA^0) \longleftarrow
H_0(\cA^1) \longleftarrow
\cdots \longleftarrow
H_0(\cA^n) \longleftarrow
0
$$
By Lemma \ref{lAgeneral} the vertical strand $\cA^m$ in the Gobelin $\cG$
is quasi-isomorphic to the complex $\Kphimd\otimes\cO_V$ which is equal to 
\[
	0 \xleftarrow{} 
	\Omega^m \otimes \cO_V \xleftarrow{\wedge df}
	\Omega^{m-1} \otimes \cO_V \xleftarrow{\wedge df}
	\cdots
\]
in our situation.  We obtain
\[
	H_0(\cA^m) = \frac{\Omega^m}{df \wedge \Omega^{m-1}} \otimes \cO_V = \Omega^m_V
\]
and an equality of complexes $\Omega_V = H_0(\cA)_{\le n}$ since both are given by
contraction with $X$. Since $\Indhom (X;V,0) := \chi(\Omega_V)$ the claim follows.
\end{proof}

\begin{prop} \label{lConstantInFamilies}
Let $X_{\lambda}$ be a holomorphic family of germs of holomorphic vector fields in $\CC^{n+k}$
with  isolated singularities and tangent to the complete intersections
$V_{\lambda} := f^{-1}(\alpha(\lambda))$, 
$$d(f-\alpha(\lambda)).X_{\lambda} = c_{\lambda}.(f-\alpha(\lambda))$$ 
with $0$ an isolated singularity for $V_0$, and the other $V_{\lambda}$ smooth. Then we have

$$\Indhom (X_0;V,0) =
\sum_{X_{\lambda}(p_{\lambda,j})=0} \Indhom (X_{\lambda}; V_{\lambda}, p_{\lambda,j}) =$$
$$=
\sum_{X_{\lambda}(p_{\lambda,j})=0} \IndGSV (X_{\lambda}; V_{\lambda},p_{\lambda,j}) =
\IndGSV (X_0; V_0,0)  $$
\end{prop}

\begin{proof}
Consider the family of cut Gobelins $\cG^{\cut}_{\lambda}$ constructed for $(X_{\lambda},c_{\lambda})$, 
choose a conveniently small ball $U$ around $0$ and 
denote the double complex of sheaves  on $U$ obtained from the Gobelins by $\cG^{\cut}_{U,\lambda}$.
The free $\cO_U$ modules in the Gobelin are independent
of $\lambda$ and the morphisms are dependent on $\lambda$. 
By Proposition \ref{pCutIndex} the Euler characteristic of 
of $\tot(\cG^{\cut}_0)$ is the homological index of $X_0$ at $0$.

Now for $\lambda \neq0$, since $V_{\lambda}$ is smooth, the complex
has only non-zero homology at degree $0$ and its dimension is equal
to $\sum_{X_{\lambda}(p_{\lambda,j})=0} \IndGSV  (X_{\lambda}; V_{\lambda}, p_{\lambda,j})$, since at smooth
points the homological and the GSV-index coincide by Koszul's Theorem.
Hence the Euler characteristic
of $\tot(\cG^{\cut}_{\lambda})$ equals the above sum, for $\lambda \neq0$. 

Now the family of holomorphic complexes of sheaves
$\tot({\cG}^{\cut}_{U,\lambda})$ on $U$
is formed by free sheaves on $U$
having cohomology sheaves  supported on $\{X_{\lambda}=0\}$ and hence
the projection of the supports to $\lambda$ is a finite map. These are the hypothesis needed in the Theorem from \cite{GG02}, and we obtain as a conclusion that the Euler characteristics of
the complexes coincide for $0$ and for small values of $\lambda$.   
\end{proof}

\begin{sloppypar}

\begin{prop} \label{lFamilyExists}
Let $V_\lambda$ be a holomorphic family of complete intersection germs defined by $f_\lambda = (f_{1,\lambda},\dots,f_{k,\lambda})^t$
in $\CC^N$, such that $V_0$ has an isolated singularity in $0$. 
Then there exists a holomorphic family of holomorphic vector fields $X_\lambda$ tangent to $V_\lambda$ such that for
$\lambda$ small enough $X_\lambda$ has isolated singularities in $\CC^N$ near $0$.
\end{prop}

\end{sloppypar}

\begin{proof}
Let $U \subset \CC^N$ be an open neighbourhood of $0$ where all components of $f_\lambda$ are convergent for small $\lambda$. From now on we denote by $V_\lambda$ and $X_\lambda$ representatives in $U$.

A vector field $X_\lambda$ is tangent to $V_\lambda$ if there exists a $k \times k$ matrix $c_\lambda$ such that
\[
	\varphi_\lambda X_\lambda = c_\lambda f_\lambda
\]
with $\varphi_\lambda$ the matrix of partial derivatives of $f$. Collecting the coefficients of $X_\lambda$ and $c_\lambda$
in one column, we obtain the matrix equation
\[
	\left(
	\begin{array}{c|ccc}
	 & f_\lambda^t & & \\
	\varphi_\lambda && \dots & \\
	&&&f_\lambda^t  
	\end{array}
	\right)
	\left(
	\begin{array}{c}
	X_\lambda \\
	-c_{\lambda,11} \\
	-c_{\lambda,12} \\
	\vdots \\
	-c_{\lambda,kk}
	\end{array}
	\right) = 0
\]
We write $(\varphi_\lambda \otimes 1 \, | \, E_k \otimes f_\lambda^t)$ for the left-hand matrix. Note that we have
\[ (A \otimes B) (C \otimes D) = (AC) \otimes (BD)
\]
for matrices $A,B,C, D$ of appropriate sizes. Using this we have 
\[
	\left(
	\varphi_\lambda \otimes 1 \left| E_k \otimes f_\lambda^t	\right) \right.
	\left(
	\begin{array}{c|c}
		E_N \otimes f_\lambda^t &  \psi_\lambda \\
	\hline
		-\varphi_\lambda \otimes E_k & 0		
	\end{array}
	\right)
	=0
\]
where $\psi_\lambda$ is the $N \times {N \choose k+1}$-matrix from the Buchsbaum-Rim complex $\cC^1(\varphi_\lambda)$ presenting the kernel of $\varphi_\lambda$. (This holds since in our situation $\varphi_\lambda$ drops rank in expected codimension.)
With
\[
	M_\lambda = \left( E_N \otimes f_\lambda^t \left|  \psi_\lambda
		\right) \right.
\]
and $w \in \CC^{Nk+{N \choose k+1}}=:W$ we obtain a family $X_{\lambda,w} := M_\lambda \cdot w$ that is tangent to all $V_\lambda$. We now prove that there exists a $w \in W$ such that $X_{0,w}$ has isolated singularities
in $U$. This is done by a dimension count. Consider the incidence variety
\[
	I = \{(a,w) \,|\, M_0(a)w =0\} \subset U \times W
\]
and the natural projections
\begin{center}\mbox{
\xymatrix{
	I \ar[r]^{q} \ar[d]^{p}& W \\
	U
	}
}
\end{center}
We have three cases. For $a \in U- V_0$ at least one component of $f_\lambda(a)$ is non zero
and therefore $M_0(a)$ has full rank $N$, the codimension of $p^{-1}(a)$ is $N = \dim(U-V_0)$ and 
\[
	\dim p^{-1}(U - V_0) = \dim W.
\]
Let $S$ be the singular set of $V_0$.  If $a\in V_0 - S$ then $\varphi_0(a)$ has full rank $k$. By Buchsbaum-Rim $\psi_0$ also presents the kernel of $\varphi_0$ and therefore
$\rank M_0(a) = \rank \psi_0(a) = N-k = \dim(V_0-S)$.
It follows that
\[
	\dim p^{-1}(V_0 - S) = \dim W.
\]
Finally for $a \in S$ we do not know anything, but since $S$ is finite we still have
\[
	\dim p^{-1}(S) \le \dim W.
\]
This proves that $\dim I = \dim W$ and therefore either $q$ is not surjective (and the generic fiber 
empty) or the generic fiber of $q$ is discrete. Now
\[
	q^{-1}(w) = \{a \in U \,|\, M_0(a)w = X_{0,w}(a) = 0 \}
\]
and consequently for generic $w$, $X_{0,w}$ has only isolated singularities in $U$. By semicontinuity this holds also for $X_{\lambda,w}$ with $\lambda$ small enough.
\end{proof}

\begin{theorem} \label{tConstant0}
Let $V$ be a holomorphic complete intersection germ defined by $f_\lambda = (f_{1,\lambda},\dots,f_{k,\lambda})^t$
in $\CC^N$ with an isolated singularity in $0$. If $X$ is a holomorphic vector field tangent to $V$ such that $X$ has an isolated singularity at 0, then
\[
		\Indhom (X; V,0) = \IndGSV (X; V,0)
\]
\end{theorem}

\begin{proof}
By \cite{G} the difference 
\[
	\Indhom (X; V,0) - \IndGSV (X; V,0)
\]
is a constant that depends on $V$ but not on $X$. By Proposition \ref{lConstantInFamilies} and Proposition \ref{lFamilyExists} there
exists an $X'$ such that this difference is zero. 
\end{proof}


\section{The small Gobelin} \label{ssmallGobelin}

In this section we will consider the spectral sequence that takes horizontal homology of the
Gobelin first.

\subsection{The Horizontal Complexes $\cB^m$ of the Gobelin}

\newcommand{\Kc}{\cK_{(s\cdot c)^t}}
\newcommand{\Kcmd}{(\Kc^{m})^*}
\newcommand{\Kcm}{\Kc^{m}}
\newcommand{\Kcnm}{\Kc^{n-m}}

Denote by $\cB^m$ the complex obtained from the $m^{th}$ row of the Gobelin.
The complexes $\cB^m$ are the tensor product of two complexes, arising from the direct sum decomposition of the middle module in the syzygy (\ref{(0.3)}). One is the part of the Buchsbaum-Eisenbud
strand of the Koszul complex associated to the map
\begin{equation*}
c \colon H\otimes \cO \longrightarrow G \otimes \cO
\end{equation*}
above the splicing map and the other is the
Koszul complex associated to the sequence $X_1,\ldots,X_N$.
In this subsection we apply the Buchsbaum-Eisenbud and Koszul Theorems
to describe the homologies of the horizontal complexes in the Gobelin.

Denote by $\cK_X$ the Koszul complex obtained from the morphism $X\colon \cO \longrightarrow F \otimes \cO$
and by $\cK^*_X$ its dual:
\begin{equation*}
0 \longleftarrow
  \cO  \longleftarrow
F^*  \otimes  \cO  \longleftarrow
\Lambda^2 F^* \otimes  \cO \longleftarrow
\cdots \longleftarrow
\Lambda^N F^* \otimes  \cO \longleftarrow
0
\end{equation*}
whose morphisms are contractions with $X^t=(X_1,\ldots,X_N)$.
Its 0-homology group is
$$H_0(\cK_X^*) = \BB := {\frac{\cO}{(X_1,\ldots,X_N)}}.$$
If $X_1, \ldots , X_N$ is a regular $\cO$-sequence, then by Koszul's Theorem
 $H_j(\cK_X^*)=0$ for $j>0$.

Consider  the Koszul complex $\Kc$ obtained from the morphism
$(c\cdot s)^t\colon \widetilde \cO \longrightarrow H^*\otimes \widetilde \cO(1)$ tensored with $\widetilde \cO(m-k)$:
\begin{equation*}
{
0 \leftarrow   
\Lambda^k H^*  \otimes \widetilde \cO(m) \leftarrow
\Lambda^{k-1} H^* \otimes \widetilde \cO(m-1)  \leftarrow
\cdots \leftarrow
\widetilde \cO(m-k)
\leftarrow 0.}
\end{equation*}
We denote the complex of its global sections
\begin{equation*}
{
0 \leftarrow  
S_{m}G \otimes \Lambda^k H^* \cO  \leftarrow
S_{m-1}G \otimes \Lambda^{k-1} H^* \otimes \cO  \leftarrow
\cdots \leftarrow
S_{m-k}G \otimes \cO
\leftarrow 0}
\end{equation*}
by $\Kcm$ 
 and  by $\Kcmd$
the dual complex
\begin{equation*}
{
0 \leftarrow  
D_{m-k}G^* \otimes \cO  \leftarrow
\cdots \leftarrow
D_{m-1}G^* \otimes \Lambda^{k-1} H \otimes \cO \leftarrow
D_{m}G^*\otimes \Lambda^k H \otimes \cO
\leftarrow 0}
\end{equation*}
For $m \ge k$ these complexes have length $k$.

\begin{lemma} \label{lBgeneral}
For all $m \ge 0$ we have
\begin{enumerate} 

 \item $\cB^m = \Kcmd \otimes_\cO\cK_X^*$ 

\item
If $X_1,\ldots,X_N$ is a regular $\cO$-sequence then the complex $\cB^m$ is quasi-isomorphic to the complex
$\Kcmd \otimes_\cO \BB $ and
$H_j(\cB^m) = 0$ for $j>\min(k,m)$.

\end{enumerate}
\end{lemma}

\begin{proof}
Same proof as for Lemma \ref{lAgeneral}.
\end{proof}


\subsection{The Small Gobelin}

Over the ring $\BB :=\cO/(X_1,\ldots,X_N)$ the identity $\varphi X = cf$ reduces to $cf=0$.
Using the notation
$\widetilde \BB:= \BB \otimes \cO_\P$ and $\widetilde \BB(1):= \BB \otimes \cO_\P(1)$
it gives rise to a smaller syzygy
\begin{equation*}
\widetilde\BB
\xrightarrow{ (s\cdot c)^t}
H^* \otimes\widetilde\BB(1)
\xrightarrow{-f^t}
\widetilde\BB(1).
\end{equation*}
Weaving the two associated Koszul complexes, taking global sections and dualizing we obtain
a smaller Gobelin, which we call $\cG_\BB$.  Its colums are 
$\cA^m_\BB = D_m G^* \otimes \cK^*_f \otimes \BB$ and its rows are 
$\cB^m_\BB = \Kcmd \otimes \BB$.

\addvspace{3mm}

\begin{theorem} \label{tB}
If $X_1,\ldots,X_N$ is an $\cO$-regular sequence, then the total Gobelin complexes
 $\tot(\cG)$ and  $\tot(\cG_\BB)$ are quasi-isomorphic.
\end{theorem}

\begin{proof}
By Part 2 of Lemma~\ref{lBgeneral} the horizontal rows $\cB^m_\BB$ and
$\cB^m$ are quasi-isomorphic. Hence the horizontal homology of $\cG$ and $\cG_\BB$
coincides. By construction the vertical maps between these homologies are also the same.
Consequently the spectral sequence of $\cG$ and of $\cG_\BB$ where we do
the horizontal homology first shows that the two total complexes are quasi-isomorphic. 
\end{proof}

\noindent Note that the small Gobelin $\cG_\BB$ is much simpler than the big Gobelin $\cG$ (see Fig.~\ref{fsGobelin}).

\begin{figure}[h]
{\tiny
$$
\hskip -15mm
\begin{array}{*{18}{c@{\,}}}

&
&&
&&
&&
&&
& \downarrow &
& \downarrow &
\\
\\

&
&&
&&
&&
&0 & \longleftarrow
& D_4G^* \otimes \BB & \longleftarrow
&
D_5G^*\otimes H\otimes  \BB 
\\ \\

&
&&
&&
&&
& \downarrow &
& \downarrow &
& \downarrow &

\\ \\

&&
&&&& 0 & \longleftarrow
&D_3G^* \otimes \BB & \longleftarrow
&
D_4G^*\otimes H\otimes  \BB
 & \longleftarrow
&
D_5G^*\otimes\Lambda^2H \otimes \BB
&

\\ \\

&
&&
&&  
& \downarrow &
& \downarrow &
& \downarrow &
& \downarrow &

\\ \\

&&&&
 0 & \longleftarrow
&D_2G^*\otimes\BB & \longleftarrow
&
D_3G^*\otimes H  \otimes \BB
 & \longleftarrow
&
D_4G^*\otimes\Lambda^2H \otimes\BB
 & \longleftarrow
& 0

\\ \\

 &&&
& \downarrow &
& \downarrow &
& \downarrow &
& \downarrow &

\\ \\

&&0          & \longleftarrow
& G^* \otimes \BB & \longleftarrow
&
D_2G^*\otimes H \otimes \BB 
 & \longleftarrow
&
D_3G^*\otimes\Lambda^2H \otimes \BB
 & \longleftarrow
& 0

\\ \\

&&\downarrow  &
& \downarrow &
& \downarrow &
& \downarrow &

\\ \\

0          & \longleftarrow &
 \BB & \longleftarrow
&
G^*\otimes  H\otimes \BB 
 & \longleftarrow
&
D_2G^*\otimes\Lambda^2H \otimes \BB
 & \longleftarrow 
&0

\\ \\

&& \downarrow &
& \downarrow &
& \downarrow &

\\ \\

0          & \longleftarrow &
 H\otimes \BB & \longleftarrow
&

G^*\otimes\Lambda^2H\otimes \BB

& \longleftarrow
& 0

\\ \\

&
& \downarrow &
& \downarrow &

\\ \\

0          & \longleftarrow &
 \Lambda^2H\otimes \BB & \longleftarrow
& 0

\\ \\

&
& \downarrow 
\\ \\

&
& 0
\\ \\
\end{array}
$$}
\caption{The lower left hand part of the small Gobelin
$\cG_\BB$ for $N=4$, $k=2$, $l=2$
beginning at $(0,0)$.} \label{fsGobelin}
\end{figure}


\begin{prop} \label{pRank}
If $k=l$ then $\rank_\BB (\tot \cG_\BB)_i = {k + i -1 \choose i}$.
\end{prop}

\begin{proof}
We have
\[
	(\tot \cG_\BB)_i 
	= \bigoplus_{j=0}^i (\cG_\BB)_{j,i-j}
	= \bigoplus_{j=0}^i  D_{j} G^* \otimes \bigwedge^{k+j-(i-j)} H 
	\cong \bigoplus_{j=0}^i D_{j} G^* \otimes \bigwedge^{i-2j} H 
\]
A basis of $D_j G^*$ is given by monomials $x^\alpha$ of degree $j$ in $k$ variables. A basis of
$\bigwedge^{i-2j} H$ is given by square free monomials $x^\beta$ of degree $i-2j$ in $k$ variables.
Let $S_i$ be the $\BB$-module spanned by all monomials  of degree $i$ in $k$ variables.
We then have a natural map of free $\BB$-modules:
\[
	\begin{matrix}
	 	\mu \colon & (\tot \cG_\BB)_i &\to &S_i \\
				      & x^\alpha \otimes x^{\beta} & \mapsto & x^{2\alpha+\beta}
	\end{matrix}
\]
For an inverse of this map let $x^\alpha$ with $\alpha = (\alpha_1, \dots, \alpha_k)$ be a monomial of degree $i=|\alpha|$ in $k$ variables. Consider the parity-function 
$$
	\sigma(\alpha_l) = \left\{ \begin{matrix} 
	0 & \text{if $\alpha_l$ is odd}\\ 
	1& \text{if $\alpha_l$ is odd}
	\end{matrix} 
	\right.
$$
Then $x^{\alpha - \sigma(\alpha)}$ has only even exponents while $x^{\sigma(\alpha)}$ is square free. The inverse of $\mu$ is therefore given by 
\[
	\begin{matrix}
	 	\mu^{-1} \colon & S_i & \to &(\tot \cG_\BB)_i \\
				      & x^{\alpha} & \mapsto & 
				      x^{\frac{1}{2}(\alpha - \sigma(\alpha))} \otimes x^{\sigma(\alpha)}.\\
	\end{matrix}
\]
It follows that the rank of $(\cG_\BB)_i$ as $\BB$-module is given by the number of monomials of degree
$i$ in $k$ variables, i.e. ${ k+i-1 \choose i}$.
\end{proof}

\section{Algebraic Formulae} \label{sFormulae}

In the geometric situation of Section \ref{sConstant0} the ring $\BB:=\cO/(X_1,\ldots,X_N)$ is
a finite dimensional $\KK$-vector space. Using the small Gobelin we can reduce the calculation of 
the homological index to linear algebra:

\begin{theorem}  \label{tSingularity}
Let $(V ,0) \subset (\CC^N,0)$ be a complete intersection of dimension $n$ defined by
a regular sequence $f = (f_1,\dots,f_k)$ of holomorphic function germs and 
$X = (X_1,\dots,X_N)$ the germ of a holomorphic vector field  tangent to $V$ in $\CC^N$.
If $V$ and $X$ have an isolated
singularity at $0$ then 
$$
H_i(\Omega_V) = \left\{ 
\begin{matrix}
H_i(\tot(\cG_\BB)) & \text{for $i \le n-1$} \\
\cO_V/I_k(\varphi) & \text{for $i=n$}
\end{matrix} \right.
$$
\end{theorem}

\begin{proof} 
By Theorem \ref{tCut} we have $H_i(\Omega_V) = H_i(\tot(\cG^{\cut}_{\le n}))$ for $i\le n$.
Since the columns $\cA^m$, $m\le n$, of $\cG$ and $\cG^{\cut}_{\le n}$ are resolutions, taking horizontal cohomology first shows that $H_i(\tot(\cG^{\cut}_{\le n})) = H_i(\tot(\cG))$ for $i \le n-1$.
Theorem \ref{tB} now gives $H_i(\tot(\cG)) = H_i(\tot(\cG_{\BB}))$. 

For $i=n$ we denote  by $T\Omega_V$ the subcomplex of torsion sheaves of $\Omega_V$ and by $\widetilde{\Omega}_V$ the quotient complex of torsion free sheaves.
By \cite{Gr75} we have $T\Omega_V^i=0$
for $i < n$. 
Since $X$ and $V$ have an 
isolated singularity at $0$ the last map of $\widetilde{\Omega}_V$ is injective and we have
$$
	H_n(\Omega_V) = T\Omega_V^n.
$$
By \cite{Gr75, Gr80} we know $T \Omega_V^n = \cO_V/I_k(\varphi)$.
\end{proof}

The dimension of $\cO_V/I_k(\varphi)$ is a numerical invariant of $V$ and is the invariant $\tau'$  of  Greuel \cite{Gr80}.

\begin{corol}  \label{cSingularity}
In the situation of Theorem \ref{tSingularity} we have
\[
\Indhom (X; V,0) = \left( \sum_{j=0}^{n-2} (-1)^j {k-1+j \choose j} \right) \dim \BB
+(-1)^{n-1} \dim \coker \gamma_{n-1} +(-1)^n  \tau'
\]
where  $\gamma_{n-1}$ is the $(n-1)$-st map of $\tot(\cG_\BB)$.
\end{corol}

\begin{proof}
Cut the total complex of the small Gobelin at the $(n-1)$-st map
to obtain a shorter complex $\tot(\cG_\BB)_{\le n-1}$. We have
on the one hand
\[
	\chi(\tot(\cG_\BB)_{\le n-1})= 
	\sum_{j=0}^{n-2} (-1)^{j} \dim H_j(\tot \cG_{\BB}) + (-1)^{n-1} \dim \ker \gamma_{n-2}
\]
and on the other hand
\[
	\chi(\tot(\cG_\BB)_{\le n-1})= 
	\sum_{j=0}^{n-2} (-1)^{j} \rank \tot(\cG_{\BB})_j + (-1)^{n-1}\dim \tot(\cG_\BB)_{n-1} .
\]
Since 
\begin{align*}
     \dim H_{n-1}(\tot(\cG_\BB))
     &= \dim \ker \gamma_{n-2} - \dim \Img \gamma_{n-1}\\
     &= \dim \ker \gamma_{n-2} - \dim \tot(\cG_\BB)_{n-1} + \dim \coker \gamma_{n-1},
\end{align*}
we obtain the desired formula by Proposition~\ref{pRank} and Theorem~\ref{tSingularity}.
\end{proof}

We will now apply this general formula in special situations.

\subsection{Codimension 1 in $\CC^N$}

Here $k=1$. The total complex of $\cG_\BB$ is then
\[
0 \xleftarrow{} 
\BB \xleftarrow{f}
\BB \xleftarrow{c}
\BB \xleftarrow{f}
\BB \xleftarrow{c}
\BB \xleftarrow{}
\cdots
\]
By Theorem~\ref{tSingularity} we have
\[
	H_{i}(\Omega_V) = 
	\left\{
		\begin{array}{cl} 			
		\BB/(f) & \text{if $i=0$} \\
		\Ann_\BB (c) / (f) & \text{if $i$ is even} \\
		\Ann_\BB (f) / (c) & \text{if $i$ is odd} \\
		\end{array}
	\right.
\]
and by Corollary~\ref{cSingularity}
\[
	\Indhom (X; V,0) = \left\{
		\begin{array}{cl} 			
		\dim \BB - \dim \BB/(c) + \tau' & \text{for $n$ even} \\
		\dim \BB/(f) - \tau' & \text{for $n$ odd}. \\
		\end{array}
	\right.
\]
We recover thus the results of \cite{G}.

		
\subsection{Codimension 2 in $\CC^N$}

Here $k=2$. The total complex of $\cG_\BB$ is in this case 
\[
0 \xleftarrow{}
\BB^1
\xleftarrow{\left(
\begin{smallmatrix}{{f}}_{1}&
{{f}}_{{2}}\\
\end{smallmatrix}
\right)}
\BB^2
\xleftarrow{\left(
\begin{smallmatrix}{-{{f}}_{{2}}}&
{{c}}_{11}&
{{c}}_{{2}1}\\
{{f}}_{1}&
{{c}}_{1{2}}&
{{c}}_{{2}{2}}\\
\end{smallmatrix}
\right)}
\BB^3
\xleftarrow{\left(
\begin{smallmatrix}{-{{c}}_{1{2}}}&
{-{{c}}_{{2}{2}}}&
{{c}}_{11}&
{{c}}_{{2}1}\\
{{f}}_{1}&
0&
{{f}}_{{2}}&
0\\
0&
{{f}}_{1}&
0&
{{f}}_{{2}}\\
\end{smallmatrix}
\right)}
\BB^4
\xleftarrow{\left(
\begin{smallmatrix}{-{{f}}_{{2}}}&
0&
{{c}}_{11}&
{{c}}_{{2}1}&
0\\
0&
{-{{f}}_{{2}}}&
0&
{{c}}_{11}&
{{c}}_{{2}1}\\
{{f}}_{1}&
0&
{{c}}_{1{2}}&
{{c}}_{{2}{2}}&
0\\
0&
{{f}}_{1}&
0&
{{c}}_{1{2}}&
{{c}}_{{2}{2}}\\
\end{smallmatrix}
\right)}
\BB^5
\xleftarrow{}
\cdots
\]
Or more generally 
\[
	\cdots 
	\xleftarrow{} 
	\BB^{2i-1} 
	\xleftarrow{\varphi_i} 
	\BB^{2i} 
	\xleftarrow{\psi_i} 
	\BB^{2i+1}
	\xleftarrow{}
	\cdots
\]
with
\renewcommand{\arraystretch}{0.5}
\[
{\scriptsize \scriptstyle
	\varphi_i = 
	\left(
	\begin{array}{*3{r@{\,}}r|*3{r@{\,}}r}
	-c_{12} & -c_{22}&&& c_{11} & c_{21}&& \\
	 &\ddots & \ddots && & \ddots & \ddots & \\
	& & -c_{12} & -c_{22}&& & c_{11} & c_{21}  \\
	&&&&&&\\\hline
	&&&&&&\\
	f_1 &&&& f_2 &\\
	& \ddots &&&& \ddots \\
	&& \ddots &&&& \ddots \\
	&&& f_1 &&&& f_2
	\end{array}
	\right), \quad
	\psi_i = 
	\left(
	\begin{array}{*2{c@{\,}}c|c@{\,}c@{}c@{\,}c}
	-f_2 & & &c_{11} & c_{21}&& \\
	&\ddots & && \ddots & \ddots & \\
	&& -f_2 & && c_{11} & c_{21}  \\
	&&&&&&\\\hline
	&&&&&&\\
	f_1 &&& -c_{12} & -c_{22}&\\
	&\ddots&& & \ddots & \ddots & \\
	&&f_1& & & -c_{12} & -c_{22}
	\end{array}
	\right).}
\]	
\renewcommand{\arraystretch}{1.0}
From Corollary~\ref{cSingularity} we obtain
\newcommand{\codimtwo}{
\[
	\Indhom (X; V,0) = 
	\left\{
		\begin{array}{cl} 			
		i(\dim \BB)-\coker \psi_i + \tau' & \text{for $n=2i$ even} \\
		(1-i)(\dim \BB)+\coker \varphi_i - \tau' & \text{for $n=2i-1$ odd.} \\
		\end{array}
	\right.
\]		
}
\codimtwo


\subsection{Complete Intersection Curves and Surfaces}

If $k \ge 2$ the total complex of $\cG_\BB$ is
\[
0
\xleftarrow{}
\BB
\xleftarrow[\gamma_0]{\left(\begin{smallmatrix}{{f}}_{1}&
      {{f}}_{{2}}&
	\cdots &
      {{f}}_{{k}}\\
      \end{smallmatrix}\right)}
\BB^k
\xleftarrow[\gamma_1]{\left(\begin{smallmatrix}{{f}}_{{2}}& 0 & f_3 &
      {\cdots}&
      0& 0 &
      {{f}}_{{k}}&
      {{c}}_{1,1}&
      {\cdots}&
      {{c}}_{{k},1}\\
      {-{{f}}_{1}}& f_3  & 0 &
      {\cdots}&
      0 & {{f}}_{{k}}&
      0&
      {{c}}_{1,{2}}&
      {\cdots}&
      {{c}}_{{k},{2}}\\
       0 &  -f_2& -f_1  &
      {\cdots}&
      f_k & 0&
      0&
      {{c}}_{1,{3}}&
      {\cdots}&
      {{c}}_{{k},{3}}\\
    \cdot & \cdot & \cdot & & \cdot & \cdot & \cdot & \cdot & & \cdot\\
      \cdot & \cdot & \cdot & & \cdot & \cdot & \cdot & \cdot & & \cdot\\
      0& 0 & 0 &
      {\cdots}&
      -f_3 & {-{{f}}_{{2}}}&
      {-{{f}}_{1}}&
      {{c}}_{1,{k}}&
      {\cdots}&
      {{c}}_{{k},{k}}\\
      \end{smallmatrix}\right)}
\BB^{{k \choose 2}+k}
\xleftarrow{}
\cdots
\]
If $V$ is a complete intersection curve then 
\[
	\Indhom (X; V,0) = \dim \BB/(f)-\tau'.
\]
This extends results of O.~Klehn \cite{Klehn}. If $V$ is a complete intersection surface we obtain
\[
	\Indhom (X; V,0) = \dim \BB - \coker (\gamma_1) + \tau'
\]
by Corollary \ref{cSingularity}.


\subsection{Examples}

\noindent We conclude with two explicit examples.
%
%
\begin{exam} 
Consider the homogeneous complete intersection surface $V \subset \CC^4$ 
defined by  $x^{2}+y^{2}+{z} {w}={x} {y}+z^{2}+w^{2}=0$, and the vector field
$$
X = \left(y{z}-z^{2}-w^{2},
      y^{2}-{x} {z}+\frac{1}{2}{y} {w},
      x^{2}+{y} {z}+\frac{3}{4} {z} {w},
      {y} {w}+\frac{1}{4} w^{2}\right)
$$
which is tangent to $V$ with tangency relation
\[
\varphi \cdot X = 
\left(\begin{smallmatrix}2 x&
2 y&
{w}&
{z}\\
{y}&
{x}&
2 z&
2 w\\
\end{smallmatrix}\right)
\cdot
\left(\begin{smallmatrix}{y} {z}-z^{2}-w^{2}\\
y^{2}-{x} {z}+\frac{1}{2} {y} {w}\\
x^{2}+{y} {z}+\frac{3}{4} {z} {w}\\
{y} {w}+\frac{1}{4}  w^{2}\\
\end{smallmatrix}\right)
=
\left(\begin{smallmatrix}2 y+{w}&
{-2 x}\\
{z}&
{y}+\frac{1}{2} w\\
\end{smallmatrix}\right)
\cdot
\left(\begin{smallmatrix}x^{2}+y^{2}+{z} {w}\\
{x} {y}+z^{2}+w^{2}\\
\end{smallmatrix}\right)
= c \cdot f.
\]
The total complex of the small Gobelin $\cG_\BB$ constructed from this relation is
\[
0 \xleftarrow{}
\BB^1
\xleftarrow{
\left(\begin{smallmatrix}x^{2}+y^{2}+{z} {w}&
{x} {y}+z^{2}+w^{2}\\
\end{smallmatrix}\right)
}
\BB^2
\xleftarrow{
\left(\begin{smallmatrix}-{x} {y}-z^{2}-w^{2}&
-2 y-w&
{-z}\\
x^{2}+y^{2}+{z} {w}&
2 x&
-y-1/2 w\\
\end{smallmatrix}\right)
}
\BB^3
\xleftarrow{}
\cdots
\]
where $\BB = \cO/(X_1,\dots,X_4)$ is a $16$-dimensional vector space.
We obtain $h_0(\Omega_V) =h_0(\cG_\BB) = 9$ and $h_1(\Omega_V) = h_1(\cG_\BB)=4$.
Since $\tau' =  \dim \cO_V/I_2(\varphi) = 7$ we have
\[
	\IndGSV (X; V,0) = \Indhom (X; V,0) = 9-4+7 = 12.
\]
(The calculations for this example were done using the computer algebra system {\sc Macaulay2}
\cite{M2}.)
\end{exam}

\begin{exam}  \label{Ex2}
Consider the complete intersection surface $V \subset \CC^4$ 
defined by  $x^{3}+y^{2}+{z} {w}={x} {y}+z^{2}+w^{2}=0$, and the vector field
$$
X =  \left(\begin{smallmatrix}768 {x} {z} {w}-192 {x} w^{2}+4 x^{2}-4 {x} {y}-16 {z} {w}+4 w^{2}\\
      24 {x} {z} {w}+1152 {y} {z} {w}-96 {x} w^{2}-288 {y} w^{2}+8 {x} {y}-6 y^{2}-2 z^{2}+2 w^{2}\\
      -12 x^{2} {w}+960 z^{2} {w}-288 {z} w^{2}-192 w^{3}+{x} {y}+7 {x} {z}-5 {y} {z}+z^{2}+w^{2}\\
      48 x^{2} {w}-48 {x} {y} {w}+1152 {z} w^{2}-288 w^{3}-{x} {y}+8 {y} {z}-z^{2}+5 {x} {w}-7 {y} {w}-w^{2}\\
      \end{smallmatrix}\right)
$$
which is tangent to $V$ with tangency relation
\begin{align*}
\varphi \cdot X &= 
\left(\begin{smallmatrix}3 x^{2}&
2 y&
{w}&
{z}\\
{y}&
{x}&
2 z&
2 w\\
\end{smallmatrix}\right)
\cdot
\left(\begin{smallmatrix}768 {x} {z} {w}-192 {x} w^{2}+4 x^{2}-4 {x} {y}-16 {z} {w}+4 w^{2}\\
24 {x} {z} {w}+1152 {y} {z} {w}-96 {x} w^{2}-288 {y} w^{2}+8 {x} {y}-6 y^{2}-2 z^{2}+2 w^{2}\\
-12 x^{2} {w}+960 z^{2} {w}-288 {z} w^{2}-192 w^{3}+{x} {y}+7 {x} {z}-5 {y} {z}+z^{2}+w^{2}\\
48 x^{2} {w}-48 {x} {y} {w}+1152 {z} w^{2}-288 w^{3}-{x} {y}+8 {y} {z}-z^{2}+5 {x} {w}-7 {y} {w}-w^{2}\\
\end{smallmatrix}\right)
\\ &=
\left(\begin{smallmatrix}2304 {z} {w}-576 w^{2}+12 x-12 y&
-192 w^{2}+4 y-z+{w}\\
0&
1920 {z} {w}-576 w^{2}+12 x-10 y+2 z-2 w\\
\end{smallmatrix}\right)
\cdot
\left(\begin{smallmatrix}x^{3}+y^{2}+{z} {w}\\
{x} {y}+z^{2}+w^{2}\\
\end{smallmatrix}\right)
= c \cdot f.
\end{align*}
Note that the equations for $V$ are not  weighted homogeneous.
The total complex of the small Gobelin $\cG_\BB$ constructed from this relation is
\[
0 \xleftarrow{}
\BB^1
\xleftarrow{
\left(\begin{smallmatrix}x^{3}+y^{2}+{z} {w}&
{x} {y}+z^{2}+w^{2}\\
\end{smallmatrix}\right)
}
\BB^2
\xleftarrow{
\left(\begin{smallmatrix} {x} {y}+z^{2}+w^{2}&
2304 {z} {w}-576 w^{2}+12 x-12 y&
0\\
-x^{3}-y^{2}-{z} {w}&
-192 w^{2}+4 y-z+{w}&
1920 {z} {w}-576 w^{2}+12 x-10 y+2 z-2 w\\
\end{smallmatrix}\right)
}
\BB^3
\xleftarrow{}
\cdots
\]
where $\BB = \cO/(X_1,\dots,X_4)$ is a $16$-dimensional vector space.
We obtain $h_0(\Omega_V) =h_0(\cG_\BB) = 9$ and $h_1(\Omega_V) = h_1(\cG_\BB)=4$.
Since $\tau' =  \dim \cO_V/I_2(\varphi) = 8$ we have
\[
	 \indV (X) = 9-4+8 = 13.
\]
The calculations for this example were done by the computer algebra system {\sc Singular} \cite{Si} using \cite{gobelinweb}. We are grateful to A.~Fr\"uhbis-Kr\"uger for her help.
\end{exam}


{\bf \small
\vskip 2mm
\noindent
Hans-Christian Graf von Bothmer (bothmer@math.uni-hannover.de)

\noindent
Wolfgang Ebeling (ebeling@math.uni-hannover.de)

\vskip 1mm

\noindent
Leibniz Universit\"at Hannover, Institut f\"ur Algebraische Geometrie, Postfach 6009, D-30060 Hannover, Germany

\vskip 3mm
\noindent
Xavier G\'omez-Mont (gmont@cimat.mx)

\noindent
CIMAT,
A.P. 402,
Guanajuato, 36000, M\'exico}

\end{document}